\documentclass[hidelinks,10pt]{elsarticle}
\usepackage[paperwidth=7in, paperheight=10.0in, margin=.675in]{geometry}
\usepackage[utf8]{inputenc}
\usepackage{graphicx}
\usepackage{float}
\usepackage{color}
\usepackage{microtype}
\usepackage{hyperref}

\usepackage{subfigure}
\usepackage{subfigmat}
 \usepackage{amssymb}
\usepackage{amsmath}
 \usepackage{amsthm}

\begin{document}

\begin{frontmatter}
\title{An upwind method for genuine weakly hyperbolic systems}

 \author[1]{Naveen Kumar Garg}
 \ead{garg.naveen70@gmail.com, navin10@math.iisc.ernet.in}
  \author[2]{Michael Junk}
  \ead{Michael.Junk@uni-konstanz.de}
 \author[3]{S.V. Raghurama Rao}
 \ead{raghu@aero.iisc.ernet.in}
 \author[4]{M. Sekhar}
 \ead{muddu@civil.iisc.ernet.in}
 \address[1]{Research Scholar, IISc Mathematical Initiative (IMI), Indian Institute of Science, Bangalore, India}
 \address[2]{Professor, Fachbereich Mathematik und Statistik, Universit\"at Konstanz, Germany}
 \address[3]{Associate Professor, Department of Aerospace Engineering, Indian Institute of Science, Bangalore, India}
 \address[4]{Professor, Department of Civil Engineering, Indian Institute of Science, Bangalore, India}
\begin{abstract}
 In this article, we attempted to develop an upwind scheme based on Flux Difference Splitting using Jordan canonical forms to simulate genuine weakly hyperbolic systems. Theory of Jordan Canonical Forms is being used to complete defective set of linear independent eigenvectors. Proposed FDS-J scheme is capable of recognizing various shocks accurately.
\end{abstract}

\end{frontmatter}


\section{Introduction}
Central and upwind discretization schemes are the popular categories of numerical methods for simulating hyperbolic conservation laws.  A system is said to be hyperbolic if its Jacobian matrix has all real eigenvalues with complete set of linearly independent eigenvectors.  Upwind schemes based on Flux Difference Splitting (FDS) are usually more accurate than others.  Two popular schemes belonging to this category are the approximate Riemann solvers of Roe \cite{Roe} and Osher \cite{Osher} and these are heavily dependent on eigenvector structure. Thus their applications are limited to systems which have complete set of linearly  independent eigenvectors.  Several other numerical schemes too are dependent strongly on eigenstructure and thus share the same difficulty.    

Recently, an attempt is made \cite{Smith_et_al} to extend Roe scheme to weakly hyperbolic systems by adding a perturbation parameter $\epsilon$ to make such systems strictly hyperbolic. In this article we try to develop an upwind method based on the concept of flux difference splitting together with Jordan forms, thus naming it as  FDS-J scheme, to simulate genuine weakly hyperbolic systems.  
We use the theory of Jordan canonical forms to complete the defective set of linearly independent (LI) eigenvectors.  Pressureless gas dynamics system, which happens to be weakly hyperbolic, is considered and it is known to produce delta shocks for density variable. Next, we consider {\em Modified Burgers' System} as given in \cite{Capdeville} and for this system too delta shocks occur exactly at same locations where normal shocks occur in the primary variables.  Similarly, other types of discontinuities, namely, $\delta^{\prime}$-shocks and  $\delta^{\prime\prime}$-shocks are  observed if we further extend modified Burgers' system as given in \cite{Shelkovich} and \cite{Joseph}.  FDS-J solver is capable of recognizing these shocks accurately.  Comparison is done with simple Local Lax-Friendrichs (LLF) \cite{LLF} method.  Contribution of generalized eigenvectors is not seen directly in the final FDS-J scheme for simulating considered genuine weakly hyperbolic systems.  It is because for each considered system, all eigenvalues are equal with arithmetic multiplicity (AM) greater than one in the resulting single Jordan block for each case.         
\section{1-D Pressureless system}  
Consider the one-dimensional pressure-less gas dynamics system
\begin{equation}\label{1D_Pressureless_system}
\frac{\partial \boldsymbol{U}}{\partial t}  + \frac{\partial \boldsymbol{F} \left( \boldsymbol{U} \right)} {\partial x} \ =  \ 0
\end{equation}
where $\boldsymbol{U}$ is the conserved variable vector and $\boldsymbol{F} \left( \boldsymbol{U} \right)$ is the flux vector defined by
\begin{equation*} 
 \boldsymbol{U} = \begin{bmatrix} 
     \rho \\[0.3em]  
		 \rho u
		\end{bmatrix} \ \mbox{and} \ 
 \boldsymbol{F}(\boldsymbol{U}) = \begin{bmatrix} 
        \rho u \\[0.3em] 
				\rho u^{2} 
				\end{bmatrix} 
\end{equation*} 
This system can also be written in quasilinear form as follows.  
\begin{equation}
\frac{\partial \boldsymbol{U}}{\partial t}  + \boldsymbol{A} \frac{\partial \boldsymbol{U}} {\partial x} \ =  \ 0 \
\end{equation}
Here $A$ is  Jacobian matrix for pressure-less system  and is given by
\begin{equation*}
 \boldsymbol{A} = \begin{bmatrix}
        \ 0  & 1   \\[0.3em]
       \ -u^{2}  &  2 u 
        \end{bmatrix} \
\end{equation*} 
Eigenvalues corresponding to Jacobian matrix $A$ are $\lambda_{1} = \lambda_{2}  =  u$ and thus algebraic multiplicity (AM) of the eigenvalues is 2, so we have to find its eigenvector space to see whether $\boldsymbol{A}$ has complete set of linearly independent eigenvectors or not. The analysis of matrix $\boldsymbol{A}$ shows that given system is weakly hyperbolic as there is no complete set of linearly independent eigenvectors, with the only eigenvector being 
      \begin{equation*}
       \boldsymbol{R}_{1}  =  \begin{bmatrix}
                  \ 1    \\[0.3em]
                  \ u    
                 \end{bmatrix} 
  \end{equation*}
Since given system doesn't have a complete set of linearly independent (LI)  eigenvectors, it will be difficult to apply any upwind scheme based on either Flux Vector Splitting (FVS) method or Flux Difference Splitting (FDS) method.  But from the theory of Jordan Canonical Forms we can still recover complete set of LI generalized eigenvectors.
\section{Jordan canonical forms and FDS for Pressureless Gas Dynamics}\label{Jordan_forms}
Every square matrix is similar to a triangular matrix with all eigenvalues on its main diagonal. A square matrix is said to be  similar to a diagonal matrix only if it has a complete set of LI eigenvectors.  But every square matrix can be made similar to a {\em Jordan} matrix. An $n\times n$ matrix $\boldsymbol{J}$ with repeated eigenvalue $\lambda$ is called a $Jordan$ matrix of order \textbf{n} if each diagonal entry in a Jordan block is $\lambda$, each entry in the {\em super diagonal} is $1$ and every other entry is zero.  Here we are providing a brief procedure to reduce a given square matrix to a Jordan matrix.  
   
\subsection{Re-visit of typical cases}
   Let $\boldsymbol{A}$ be $n \times n$ matrix with $n$ real eigenvalues $\lambda_{1}, \lambda_{2}, \lambda_{3}, \cdots, \lambda_{n}$.  Now the following  typical cases may arise:
   
   Case $1$:  When all $\lambda_{i}$, where $1\leq i \leq n$, are distinct. In this case matrix $\boldsymbol{A}$ will have a complete set of LI eigenvectors and hence will be similar to a diagonal matrix.
  
  Case $2$: When some $\lambda_{i}$ are equal, {\em i.e.}, let $\lambda_{1}  \ = \  \lambda_{2} \ = \ \lambda_{3} \ = \ \cdots \ = \ \lambda_{p} \ = \ \lambda$, where $p$ is a natural number $\leq n$, any of the  following sub-cases may happen.
  
  Sub-case $1$: If algebraic multiplicity (AM) of an eigenvalue $\lambda$, which is $p$ in assumed case, is  equal to geometric multiplicity (GM), and moreover if this is true for all subsets of equal eigenvalues then square matrix $\boldsymbol{A}$ will again be similar to a unique diagonal matrix. 
  
  Sub-case $2$: Now, consider the case in which GM is strictly less than AM, in that case the LI set of eigenvectors will not be a complete one.  Here, we can pull in the theory of Jordan canonical forms to recover full LI set of generalized eigenvectors and to make given square matrix similar to a Jordan matrix which is not much different from a diagonal matrix.\\
\textbf{Definition:} A $n \times n$ matrix is called defective matrix if it doesn't possess full set of linearly independent eigenvectors. \\ 
\textbf{Procedure to find generalized eigenvectors:} In this article we mainly focus on systems which belong to the category as discussed in Sub-case $2$. If all eigenvalues of a given defective matrix are equal and further if there is only a  single Jordan block corresponding to given matrix, then following steps need to be followed to recover full set of LI generalized eigenvectors: \\
  $(i)$ For an eigenvalue $\lambda$, compute the ranks of the matrices $\boldsymbol{A}-\lambda \boldsymbol{I}$, \  \ $(\boldsymbol{A}-\lambda \boldsymbol{I})^2$, $\cdots$, and find the least positive integer $s$ such that $rank(\boldsymbol{A}-\lambda \boldsymbol{I})^s \ =  \  rank(\boldsymbol{A}-\lambda \boldsymbol{I})^{s+1}$. There will be a single Jordan block only if $s$ comes out equal to dimension of given matrix.   \\ 
$(ii)$  Once $s$ is equal to dimension of defective matrix, generalized eigenvectors can be computed from the system of equations $\boldsymbol{A}\boldsymbol{P} = \boldsymbol{P}\boldsymbol{J}$,  
where 
\begin{equation*}
 {\boldsymbol J(\lambda)} = \begin{bmatrix}
        \   \lambda   &   1    &        \\[0.3em]
       \    &     \ddots   &   \ddots  &  \\[0.3em]
       \  &    &   \ddots   &  1   \\[0.3em]
       \  &       &        &  \lambda
        \end{bmatrix}_{s\times s} \
\end{equation*}  \\
Let $\boldsymbol{P}$ equal to $[\boldsymbol{X_{1}}, \boldsymbol{X_{2}}, \boldsymbol{X_{3}}, ..., \boldsymbol{X_{s}}]$ be a set of column vectors which need to evaluated. Then,
\begin{equation}
 \boldsymbol{A}\big[\boldsymbol{X}_1, \boldsymbol{X}_2, \boldsymbol{X}_3,........., \boldsymbol{X}_{s}\big] \ = \ \big[\boldsymbol{X}_1, \boldsymbol{X}_2, \boldsymbol{X}_3,........., \boldsymbol{X}_{s}\big] \begin{bmatrix}
        \   \lambda   &   1    &        \\[0.3em]
       \    &     \ddots   &   \ddots  &  \\[0.3em]
       \  &    &   \ddots   &  1   \\[0.3em]
       \  &       &        &  \lambda
        \end{bmatrix}_{s\times s} \
\end{equation}
gives
 \begin{align}\label{formula_to_find_gen_eig}
 \begin{split}
  \boldsymbol{A}\boldsymbol{X}_{1}  \ &= \  \lambda \boldsymbol{X}_{1} \\
  \boldsymbol{A}\boldsymbol{X}_{2}  \ &= \  \lambda \boldsymbol{X}_{2} \ + \ \boldsymbol{X}_{1}  \\
  \boldsymbol{A}\boldsymbol{X}_{3}  \ &= \  \lambda \boldsymbol{X}_{3} \ + \ \boldsymbol{X}_{2}  \\
                                 \vdots     \\
\boldsymbol{A}\boldsymbol{X}_{s}  \ &= \  \lambda \boldsymbol{X}_{s} \ + \ \boldsymbol{X}_{s-1}
 \end{split}
\end{align}
Now we can compute all true and generalized eigenvectors from system of relations (\ref{formula_to_find_gen_eig}).  For present case, $u$ is repeated eigenvalue with arithmetic multiplicity (AM) of 2 and on computing the ranks of matrices $\boldsymbol{A}- u\boldsymbol{I}$, \  \ $(\boldsymbol{A}- u\boldsymbol{I})^2$ and $(\boldsymbol{A}- u\boldsymbol{I})^{3}$, we find $rank(\boldsymbol{A}- u \boldsymbol{I})^2  \ =  \ 0  \ =  \  rank(\boldsymbol{A}- u \boldsymbol{I})^{3}$. Thus $s$ will be $2$ in this case, so there will be one $Jordan$ block of order $2$.  On expanding relation $\boldsymbol{A}\boldsymbol{P} = \boldsymbol{P}\boldsymbol{J}$, we get 
\begin{equation}\label{relation_to_find_gen_eigenvectors}
 \boldsymbol{A}[\boldsymbol{X}_{1} \ \  \boldsymbol{X}_{2}] \ = \ [\boldsymbol{X}_{1} \ \ \boldsymbol{X}_{2}] [\boldsymbol{J}_{1} \ \ \boldsymbol{J}_{2}]
\end{equation}
where $\boldsymbol{X}_{i}^{\prime{s}}$ are linearly independent, $2 \times 1$, column vectors. Similarly, $\boldsymbol{J}_{i}^{\prime{s}}$ are column vectors which form Jordan matrix $\boldsymbol{J}$ and are given as 
\begin{equation}
       \boldsymbol{J}_{1}  =  \begin{bmatrix}
                  \ \lambda    \\[0.3em]
                  \ 0    
                 \end{bmatrix}_{2 \times 1}, \
       \boldsymbol{J}_{2}  =   \begin{bmatrix}
                  \ 1   \\[0.3em]
                  \ \lambda
                  
              \end{bmatrix}_{2 \times 1} \
 \end{equation}
On solving (\ref{relation_to_find_gen_eigenvectors}), we get following relations to find all eigenvectors, {\em i.e.},
\begin{align}\label{def_gen_eigenvectors_2}
 \begin{split}
  \boldsymbol{A}\boldsymbol{X}_{1}  \ &= \  \lambda \boldsymbol{X}_{1} \\
  \boldsymbol{A}\boldsymbol{X}_{2}  \ &= \  \lambda \boldsymbol{X}_{2} \ + \ \boldsymbol{X}_{1} 
 \end{split}
\end{align}
First relation of (\ref{relation_to_find_gen_eigenvectors}) gives $\boldsymbol{X}_{1} = \boldsymbol{R}_{1}$ and on using this value in second relation of (\ref{relation_to_find_gen_eigenvectors}), we get 
\begin{equation*}
  \boldsymbol{X}_{2} \ = \ \boldsymbol{R}_{2} =  \begin{bmatrix}
         \ x_{1} \\[0.3em]
         \ 1 + ux_{1} 
       \end{bmatrix}
 \end{equation*}
which will be a generalized eigenvector of the pressureless gas dynamics system and $x_{1} \in {\rm I\!R}$.         
\subsection{Formulation of a FDS scheme for Pressureless System}
System (\ref{1D_Pressureless_system}) can be written in quasi-linear form  as 
\begin{equation}\label{Quasi_form_pressureless_system}
\frac{\partial \boldsymbol{U}}{\partial t}  + \boldsymbol{A} \frac{\partial \boldsymbol{U}} {\partial x} \ =  \ 0 
\end{equation}
Now, because of the non-linearity of Jacobian matrix $\boldsymbol{A}$, it is difficult to solve above system. But locally, inside each cell, $\boldsymbol{A}$ can be made linearized to form a constant matrix $\boldsymbol{\bar{A}}$, which is now a function of left and right state variables $\boldsymbol{U}_L$ and $\boldsymbol{U}_R$, {\em i.e.}, $\boldsymbol{\bar{A}}  \ = \  \boldsymbol{\bar{A}}\left(\boldsymbol{U}_L,\boldsymbol{U}_R\right)$.  
So, (\ref{Quasi_form_pressureless_system}) becomes
\begin{equation}\label{Quasi-linearised_eqn}
\frac{\partial \boldsymbol{U}}{\partial t}  + \boldsymbol{\bar{A}} \frac{\partial \boldsymbol{U}} {\partial x} \ =  \ 0 
\end{equation}
On comparing (\ref{1D_Pressureless_system}) and (\ref{Quasi-linearised_eqn}), we get 
\begin{equation}
 d\boldsymbol{F} = \boldsymbol{\bar{A}}d\boldsymbol{U}
\end{equation}
The finite difference analogue of the above differential relation is,
\begin{equation}\label{conservation1_for_pressureless_system}
 \bigtriangleup{\boldsymbol{F}}  =  \boldsymbol{\bar{A}}\bigtriangleup{\boldsymbol{U}}
\end{equation}
where,
\begin{align}
 \begin{split}
  \bigtriangleup{\boldsymbol{F}} \ &= \ \boldsymbol{F}_R - \boldsymbol{F}_L   \\
   \bigtriangleup{\boldsymbol{U}} \ &= \ \boldsymbol{U}_R - \boldsymbol{U}_L
\end{split}
\end{align}
In the above equations, subscripts $R$ and $L$ represent the right and left states respectively. Relation (\ref{conservation1_for_pressureless_system}) ensures the conservation property. As already explained the present system is weakly hyperbolic, but on the basis of above mentioned procedure, we can construct a basis of true and generalized eigenvectors for column vector $\bigtriangleup{\boldsymbol{U}}$, i.e.,
\begin{equation}
 \bigtriangleup{\boldsymbol{U}}   \ = \  \sum_{i = 1}^{2} \bar\alpha_{i}\boldsymbol{\bar{R}}_{i} 
\end{equation}
where, $\bar{\alpha}_{i}'s$ are coefficients attached with both LI eigenvectors corresponding to given system. On using above equation in (\ref{conservation1_for_pressureless_system}), we get
\begin{equation}\label{conservation2_for_pressureless_system}
 \bigtriangleup{\boldsymbol{F}}  =  \boldsymbol{\bar{A}}\sum_{i = 1}^{2} \bar\alpha_{i}\boldsymbol{\bar{R}}_{i}
 \end{equation}
For weakly hyperbolic systems, $\boldsymbol{\bar{A}}$ is non-diagonalizable, resulting in 
\begin{equation*}
 \boldsymbol{\bar{A}} \boldsymbol{\bar{R}}_{i}  \ \neq \  \bar \lambda_{i} \boldsymbol{\bar{R}}_{i}
 \ \textrm{for some i's} \
 \end{equation*}  
We now have $\boldsymbol{\bar{R}}_{2}$ as a generalized eigenvector and 
\begin{equation*}
 \boldsymbol{\bar{A}} \boldsymbol{\bar{R}}_{1}  \ = \  \bar \lambda_{1} \boldsymbol{\bar{R}}_{1} \
 ~~\mbox{and}~~
 \boldsymbol{\bar{A}} \boldsymbol{\bar{R}}_{2}  \ = \  \bar \lambda_{2} \boldsymbol{\bar{R}}_{2}  \ + \  \boldsymbol{\bar{R}}_{1}
\end{equation*}
 On using above relations in (\ref{conservation2_for_pressureless_system}), we get
\begin{equation*} 
\bigtriangleup{\boldsymbol{F}}  \ = \   \bar\alpha_{1} \bar \lambda_{1} \boldsymbol{\bar{R}}_{1}  \ + \ \bar\alpha_{2} \bar \lambda_{2} \boldsymbol{\bar{R}}_{2}    \ + \ \bar\alpha_{2}  \boldsymbol{\bar{R}}_{1}
\end{equation*} 
We now define the standard Courant splitting for the eigenvalues as 
\begin{equation*}
 \bar \lambda^{+}_{i} - \bar \lambda^{-}_{i} = |\bar \lambda_{i}| \
 \end{equation*}
After splitting each of the eigenvalues into a positive and a negative part, $\bigtriangleup{\boldsymbol{F}}^{+}$ and $\bigtriangleup{\boldsymbol{F}}^{-}$ can be written as
\begin{equation} \label{delta_F_positive}
\bigtriangleup{\boldsymbol{F}}^{+}   \ = \   \bar\alpha_{1} \bar \lambda^{+}_{1} \boldsymbol{\bar{R}}_{1}  \ + \ \bar\alpha_{2} \bar \lambda^{+}_{2} \boldsymbol{\bar{R}}_{2}  \ + \ \bar\alpha_{2}  \boldsymbol{\bar{R}}_{1}
\end{equation}
and
\begin{equation} \label{delta_F_negative}
\bigtriangleup{\boldsymbol{F}}^{-}   \ = \   \bar\alpha_{1} \bar \lambda^{-}_{1} \boldsymbol{\bar{R}}_{1}  \ + \ \bar\alpha_{2} \bar \lambda^{-}_{2} \boldsymbol{\bar{R}}_{2}  \ + \ \bar\alpha_{2}  \boldsymbol{\bar{R}}_{1}
\end{equation}
Taking a cue from the traditional flux difference splitting methods, we now write the interface flux as     
\begin{equation}
\\ \boldsymbol{F}_{I} = \frac{1}{2} \left[\boldsymbol{F}_{L} + \boldsymbol{F}_{R} \right] - \frac{1}{2} \left[ \left(\bigtriangleup{\boldsymbol{F}}^{+} \ - \ \bigtriangleup{\boldsymbol{F}}^{-}\right) \right] 
 \end{equation}  
On using (\ref{delta_F_positive}) and (\ref{delta_F_negative}) in the upwinding part of FDS formulation for pressureless system, we get
\begin{equation}
  \ \ \bigtriangleup{\boldsymbol{F}}^{+} - \bigtriangleup{\boldsymbol{F}}^{-}  \ = \  \sum_{i = 1}^{2} \bar\alpha_{i} |\bar \lambda_{i}| \boldsymbol{\bar{R}}_{i}  
\end{equation}
Since both eigenvalues are the same, above relation becomes
\begin{equation}\label{flux_differencing}
\bigtriangleup{\boldsymbol{F}}^{+} - \bigtriangleup{\boldsymbol{F}}^{-}  \ = \   |\bar \lambda|\bigtriangleup{\boldsymbol{U}}  
\end{equation}
Now $\bigtriangleup{U_{2}}$ is equal to $\bigtriangleup({\rho u})$, which can be further expressed as
\begin{equation}
 \bigtriangleup(\rho u) \ = \ \bar{u}\bigtriangleup{\rho} \ + \  \bar{\rho}\bigtriangleup{u}
 \end{equation}
where $\bar{u}$ is some average of $u_L$ and $u_R$, $\bar{\rho}$ is another average of $\rho_L$ and $\rho_R$, both to be determined. We now have  
\begin{equation}\label{U_2_average_relation}
 \rho_{R}u_{R} - \rho_{L}u_{L}  \ = \ \bar{u}(\rho_{R} - \rho_{L}) \ + \  \bar{\rho}(u_{R} - u_{L})
 \end{equation}
We need to find average values for both density and velocity variables and both of which should satisfy relation (\ref{U_2_average_relation}) to get some meaningful solutions for interface fluxes inside each cell. 
Again consider relation $\bigtriangleup{\boldsymbol{F}} = \boldsymbol{\bar{A}}\bigtriangleup{\boldsymbol{U}}$, which in expanded form can be written as
\begin{equation} \label{conservation_eq}
  \begin{bmatrix}
        \bigtriangleup (\rho u) \\[0.3em]
       \bigtriangleup  (\rho u^{2}) 
         \end{bmatrix} \ = \
    \begin{bmatrix}
         0   & 1   \\[0.3em]
         -{\bar u}^{2} & 2 \bar{u}
        \end{bmatrix} 
        \begin{bmatrix}
         \bigtriangleup (\rho)  \\[0.3em]
       \bigtriangleup  (\rho u) 
        \end{bmatrix}
       \end{equation}
First relation is automatically satisfied for any average values.  From the second relation, we get
\begin{equation}\label{momentum_eqn}
\bigtriangleup(\rho u^{2})  \ = \   -{\bar u}^{2} \bigtriangleup (\rho)   +    2 \bar{u} \bigtriangleup  (\rho u)
\end{equation}
where 
  \begin{equation}
   \bigtriangleup (\rho)  \ = \ (\rho)_R  \ - \  (\rho)_L 
  \end{equation}
\begin{equation}
 \bigtriangleup (\rho u)  \ = \ (\rho u)_R  \ - \  (\rho u)_L
\end{equation}
\begin{equation}
 \bigtriangleup (\rho u^2)  \ = \ (\rho u^2)_R  \ - \  (\rho u^2)_L
\end{equation}
After rearrangement of terms we obtain 
\begin{equation*}
   \\   {\bar{u}}^2 \bigtriangleup (\rho) \ - \ 2 \bar{u} \bigtriangleup  (\rho u) \ + \ \bigtriangleup (\rho u^{2})\ = \ 0
  \end{equation*}
which is a quadratic equation in $\bar{u}$ the solution of which, after a little algebra, is obtained as 
\begin{equation} \label{u_bar}
\bar{u}    \ = \  \frac{\sqrt{\rho_L} u_L  \ \pm \ \sqrt{\rho_R} u_R }{\sqrt{\rho_L} \ \pm \ \sqrt{\rho_R}}
\end{equation}
We neglect the root having negative signs in both numerator and denominator as it is not physical and may become infinity as $\sqrt{\rho_R} \longmapsto \sqrt{\rho_L}$ or vice-versa. Thus average value of $u$ is defined as 
\begin{equation}
\bar{u}    \ = \  \frac{\sqrt{\rho_L} u_L  \ + \ \sqrt{\rho_R} u_R }{\sqrt{\rho_L} \ + \ \sqrt{\rho_R}}
\end{equation}
On using $\bar{u}$ in the relation (\ref{U_2_average_relation}) we get 
\begin{equation}
 \rho_{R}u_{R} - \rho_{L}u_{L}  \ = \frac{\sqrt{\rho_L} u_L  \ + \ \sqrt{\rho_R} u_R }{\sqrt{\rho_L} \ + \ \sqrt{\rho_R}}(\rho_{R} - \rho_{L}) \ + \  \bar{\rho}(u_{R} - u_{L})
 \end{equation}
Now we use $(\rho_{R} - \rho_{L}) \ = \ (\sqrt{\rho_R} + \sqrt{\rho_L})(\sqrt{\rho_R} - \sqrt{\rho_L})$ in the above equation and after rearrangement of terms, we get 
\begin{equation}  
\bar{\rho} = (\sqrt{\rho_R}\sqrt{\rho_L}) 
\end{equation} 
Since density is always positive, the average value  $\bar{\rho}$ becomes equal to $(\sqrt{\rho_R\rho_L})$. One can check that the relation (\ref{U_2_average_relation}) becomes an equation for above defined averages for both density and velocity variables.  As the interface flux is now completely defined, the final update formula in the finite volume framework is written as follows.   
 \begin{equation}
 \boldsymbol{U}^{n+1}_{j} = \boldsymbol{U}^{n}_{j} - \frac{\Delta t}{\Delta x} 
 \left[ \boldsymbol{F}^{n}_{j+\frac{1}{2}} - \boldsymbol{F}^{n}_{j-\frac{1}{2}} \right] 
 \end{equation}  
\subsection{Numerical examples}
\begin{figure}[!ht]
\centerline{%
\subfigure[]{%
\includegraphics[trim=0 5 35 5, clip, width=0.55\textwidth]{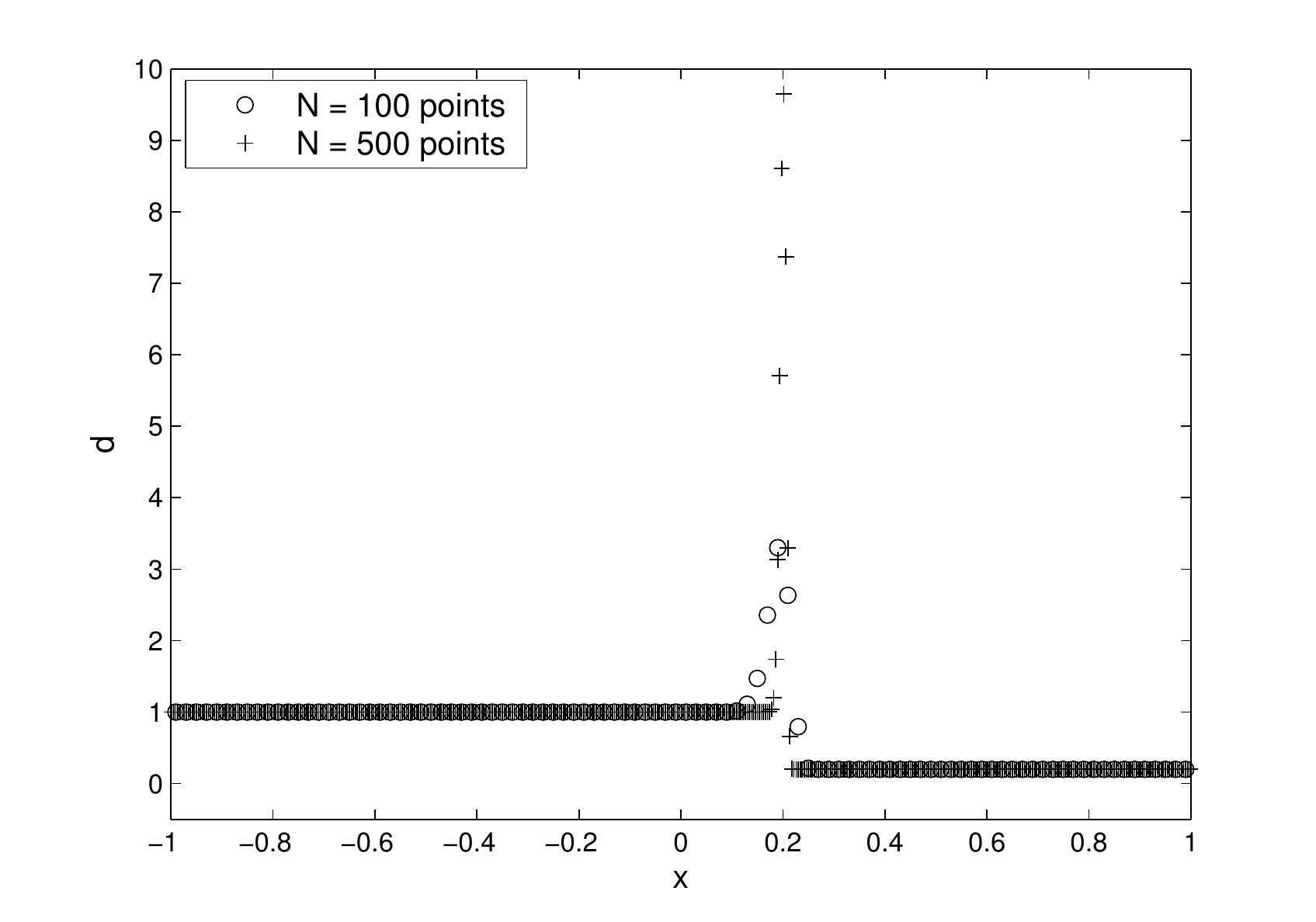}
\label{Pressureless_delta_shocks_d}
}%
\subfigure[]{%
\includegraphics[trim=0 5 35 5, clip, width=0.55\textwidth]{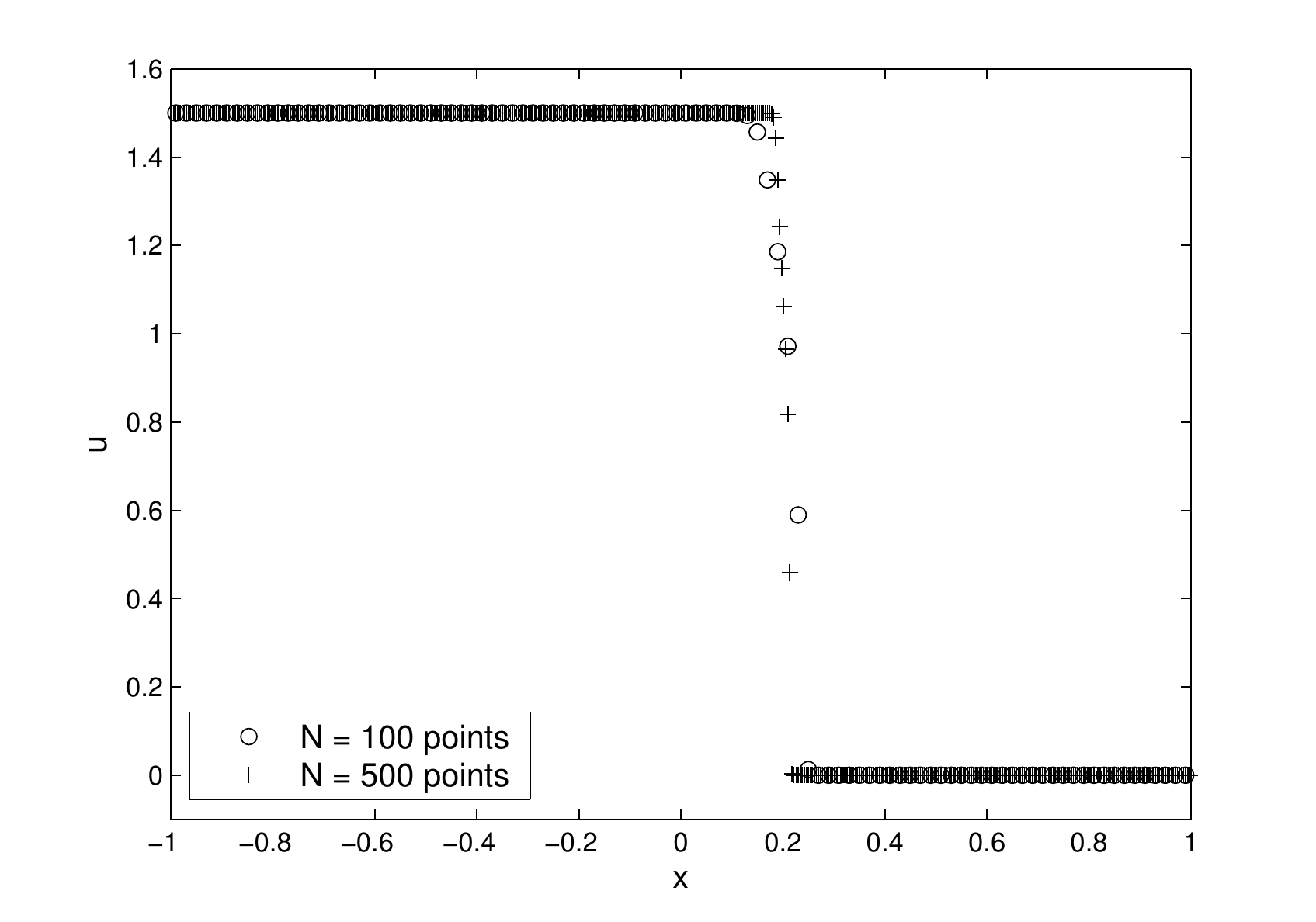}%
\label{Pressureless_delta_shocks_u}
}%
}%
\caption{ (a) Results of FDS-J scheme for pressureless system, formation of $\delta$ shocks in density variable and (b) represents formation of step discontinuity in velocity variable.}
\end{figure}

\begin{figure}[!ht]
\begin{center}
\includegraphics[trim=5 5 35 5, clip, width=0.6\textwidth]{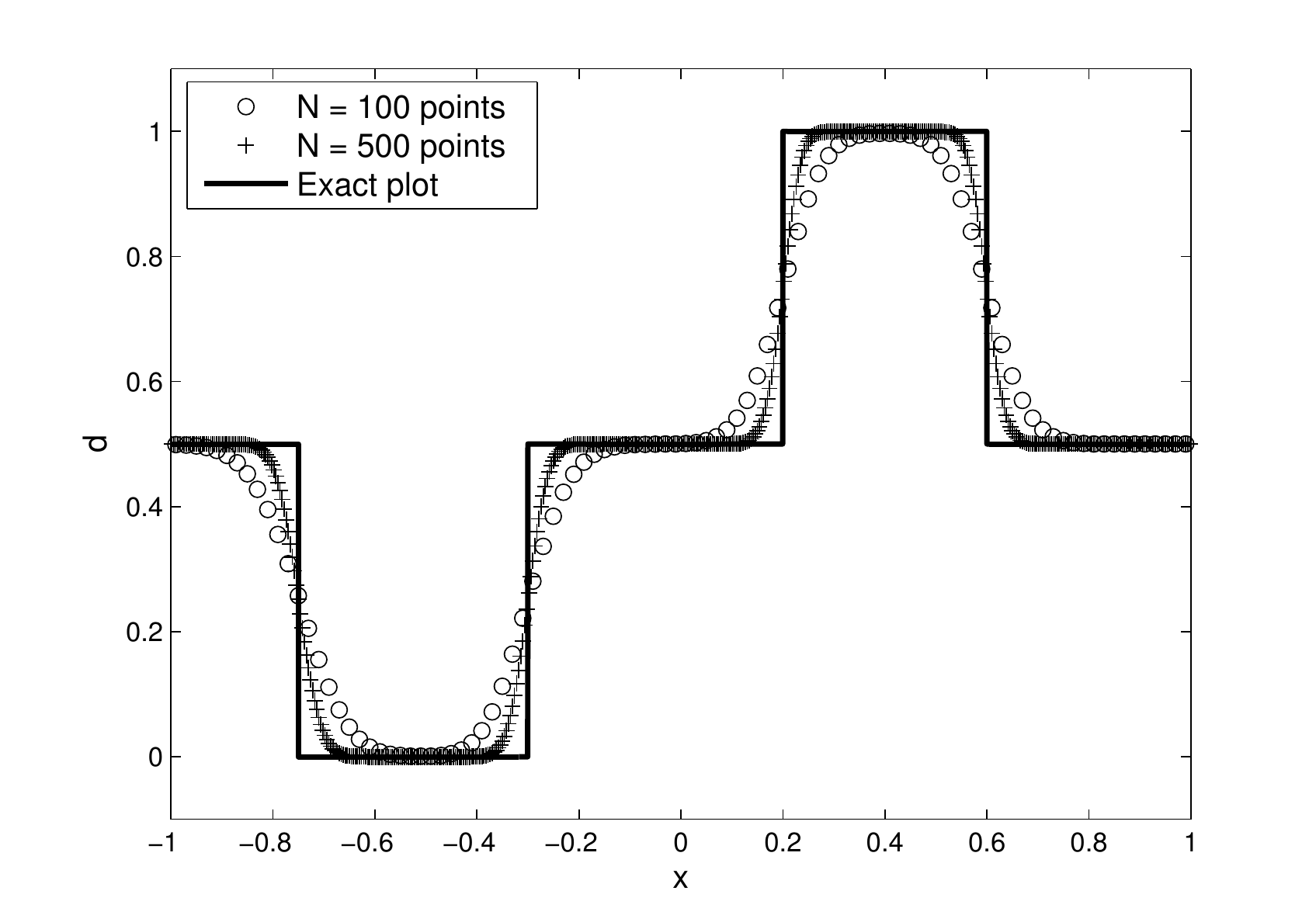}
\caption{Results of density variable for positivity problem using FDS-J scheme for pressureless gas dynamics system. }
\label{positivity_test}
\end{center}
\end{figure}

Here we consider two test cases for 1D-pressureless gas dynamics. First test case we take from \cite{Chen_&_Liu} with initial conditions being given as $(\rho_L, u_L) = (1.0, 1.5)$, $(\rho_R, u_R) = (0.2, 0.0)$ with $x_{o}=0.0$ and all solutions are obtained at final time $t=0.2$ units. In this case, a $\delta$-shock develops in density variable and our FDS-J scheme captures this feature accurately, as seen in Figure \ref{Pressureless_delta_shocks_d}. The formation of step discontinuity in velocity variable is shown in Figure 
\ref{Pressureless_delta_shocks_u}. Second test case is taken from \cite{Bouchut_Jin_&_Li}.  This test case is designed to check positivity property and maximum principle for density and velocity variables respectively.  For this problem, FDS-J scheme generates insufficient numerical diffusion.  To get meaningful solution, we use Harten's entropy fix \cite{Harten_entropy_fix} which usually increase diffusion in the scheme, {\em i.e.}, 
\begin{align}
 \begin{split}
  |\tilde{\lambda}|  \  &= \   |\lambda| \  \  \ \textrm{if}  \  \   \ |\lambda| \geq \epsilon \  \textrm{and} \ \\
  |\tilde{\lambda}|  \  &= \   \dfrac{1}{2}\big(\frac{\lambda^{2}}{\epsilon} + \epsilon \big) \ \ \textrm{if} \ \ |\lambda| < \epsilon
 \end{split}
\end{align}
for some small value of $\epsilon$.   
The  density variable plot is shown in Figure \ref{positivity_test}. 
\section{Modified Burgers' system}
Next we consider modified Burgers' system which is formed augmenting the inviscid Burgers equation with an equation obtained by taking its derivative, forming a $2\times2$ system.  Let us consider one-dimensional inviscid Burgers' equation
 \begin{equation}\label{Burgers'_1}
  u_{t}  \ + \ f_{x}(u) \ = \ 0
 \end{equation}
where, $u$ is the conserved variable and $f(u)$ is the flux function which is given by $f(u) = \frac{1}{2} u^{2}$.  On differentiating above equation $w.r.t.$ $x$, we obtain  
\begin{equation}
 (u_{t})_{x}  \ + \ (f_{x}(u))_{x} \ = \ 0
\end{equation} It further can be written as 
\begin{equation}
 (u_{x})_{t}  \ + \ (f^{\prime}(u)u_{x})_{x} \ = \ 0
 \end{equation} or
 \begin{equation}\label{Burgers'_2}
 v_{t}  \ + \ g_{x}(u) \ = \ 0
 \end{equation} where we define $v = u_{x}$ and $g(u)  = f^{\prime}(u)v $. (\ref{Burgers'_1}) and (\ref{Burgers'_2}) together form $2\times2$ system 
\begin{equation}
\frac{\partial \boldsymbol{U}}{\partial t}  \ + \ \boldsymbol{A} \frac{\partial \boldsymbol{U}}{\partial x} \ = \ 0
\end{equation}where $\boldsymbol{U}$ is column vector and $\boldsymbol{A}$ is $2\times2$ matrix, {\em i.e.},
\begin{equation*}
  \boldsymbol{U}   =  \begin{bmatrix}
         \ u \\[0.3em]
         \ v 
       \end{bmatrix}
      \ \textrm{and} \
 \boldsymbol{A} = \begin{bmatrix}
        \ u  && 0   \\[0.3em]
       \ v  &&   u 
        \end{bmatrix}
 \end{equation*}
Eigenvalues corresponding to Jacobian matrix $\boldsymbol{A}$ are $\lambda_{1} =  u = \lambda_{2}$ and thus algebraic multiplicity (AM) of the eigenvalue $u$ is 2. For $v \neq 0$, analysis of matrix $\boldsymbol{A}$ shows that given system is weakly hyperbolic as the given system has only one LI eigenvector, which is given by 
\begin{equation}
       \boldsymbol{R}_{1}  =  \begin{bmatrix}
                  \ 0    \\[0.3em]
                  \ 1    
                 \end{bmatrix}
\end{equation}
We find that there is one Jordan block of order two as $rank(\boldsymbol{A} - u \boldsymbol{I})^{2}  \ = \ 0 \ = \  rank(\boldsymbol{A} - u \boldsymbol{I})^{3}$. Like in the previous case, in order to find a generalized eigenvector we need to solve relation $\boldsymbol{A}\boldsymbol{P} = \boldsymbol{P}\boldsymbol{J}$. After a little algebra, $\boldsymbol{R}_{2}$ comes out as
\begin{equation}
 \boldsymbol{R}_{2}  =  \begin{bmatrix}
                  \ \dfrac{1}{v}    \\[0.8em]
                  \ x_{2}    
                 \end{bmatrix}
\end{equation}
where $x_{2} \in {\rm I\!R}$. 
\subsection{Formulation of FDS scheme for Modified Burgers' system}
Similar analysis like that in the pressureless gas dynamics system is valid for modified Burgers' system till equation (\ref{flux_differencing}) which is 
\begin{equation}\label{flux_differencing_2}
\bigtriangleup{\boldsymbol{F}}^{+} - \bigtriangleup{\boldsymbol{F}}^{-}  \ = \   |\bar \lambda|\bigtriangleup{\boldsymbol{U}}  
\end{equation}
In this case $\bigtriangleup{\boldsymbol{U}}$ is defined as 
\begin{equation}
  \bigtriangleup{\boldsymbol{U}} \  =  \begin{bmatrix}
         \ \bigtriangleup{u} \\[0.3em]
         \ \bigtriangleup{v} 
       \end{bmatrix}
\end{equation}
and $\bar{\lambda} = \bar{u}$. In order to solve (\ref{flux_differencing_2}) fully, we need to find average value of $u$ from relation $\bigtriangleup{\boldsymbol{F}} \ = \ \boldsymbol{\bar{A}}  \bigtriangleup{\boldsymbol{U}}$. In expanded form it can be written as  
\begin{equation} 
  \begin{bmatrix}
        \bigtriangleup (\frac{1}{2} u^{2}) \\[0.3em]
       \bigtriangleup  (u v) 
         \end{bmatrix} \ = \
    \begin{bmatrix}
         \bar u   && 0   \\[0.3em]
         \bar v &&  \bar u
        \end{bmatrix} 
        \begin{bmatrix}
         \bigtriangleup (u)  \\[0.3em]
       \bigtriangleup  (v) 
        \end{bmatrix}
       \end{equation}
From the first equation, we get
 \begin{equation}\label{1_eqn_aug_Burgers'}
  \bigtriangleup  \left(\frac{1}{2} u^{2}\right)  \ = \   \bar u  \bigtriangleup (u) 
  \end{equation}
or
\begin{equation}
 \left(\frac{1}{2} u^{2}_{R} - \frac{1}{2} u^{2}_{L}\right) \ = \  \bar{u}(u_{R} - u_{L})
\end{equation} 
if $u_{L} \neq u_{R}$, then $\bar{u} = \dfrac{(u_{L} + u_{R})}{2}$. Otherwise also $\bar{u} = \dfrac{(u_{L} + u_{R})}{2}$ in a limiting sense. Second expression ($\bar{v}$) need not be solved as interface flux requires only $\bar{u}$ to be evaluated. It is important to note that even if $v = 0$, relation (\ref{flux_differencing_2}) still holds. 
\subsection{Numerical examples}
\begin{figure}[!ht]
\centerline{%
\subfigure[]{%
\includegraphics[trim=0 5 35 5, clip, width=0.55\textwidth]{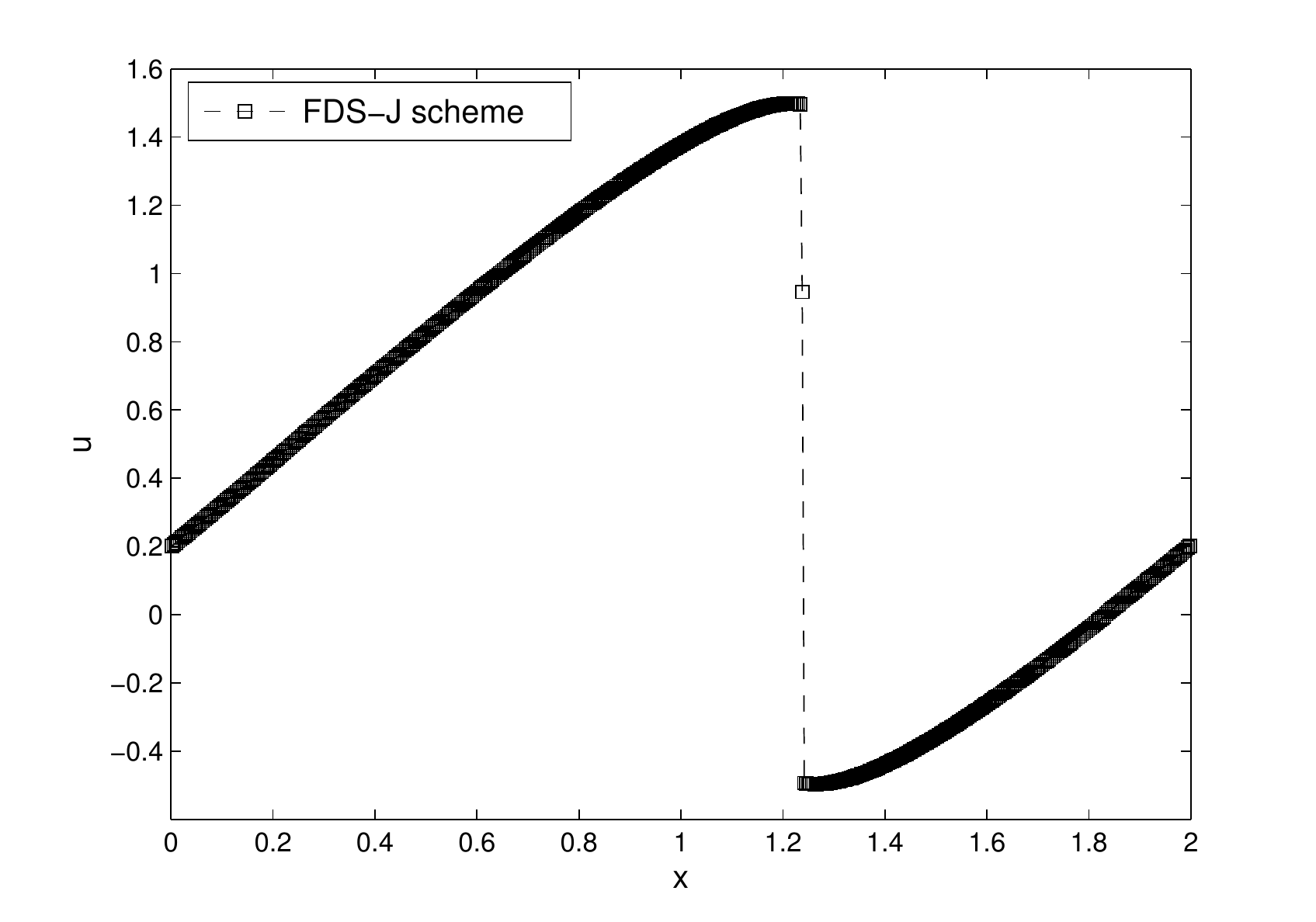}
\label{Modified Burgers'_normal_shock}
}%
\subfigure[]{%
\includegraphics[trim=0 5 35 5, clip, width=0.55\textwidth]{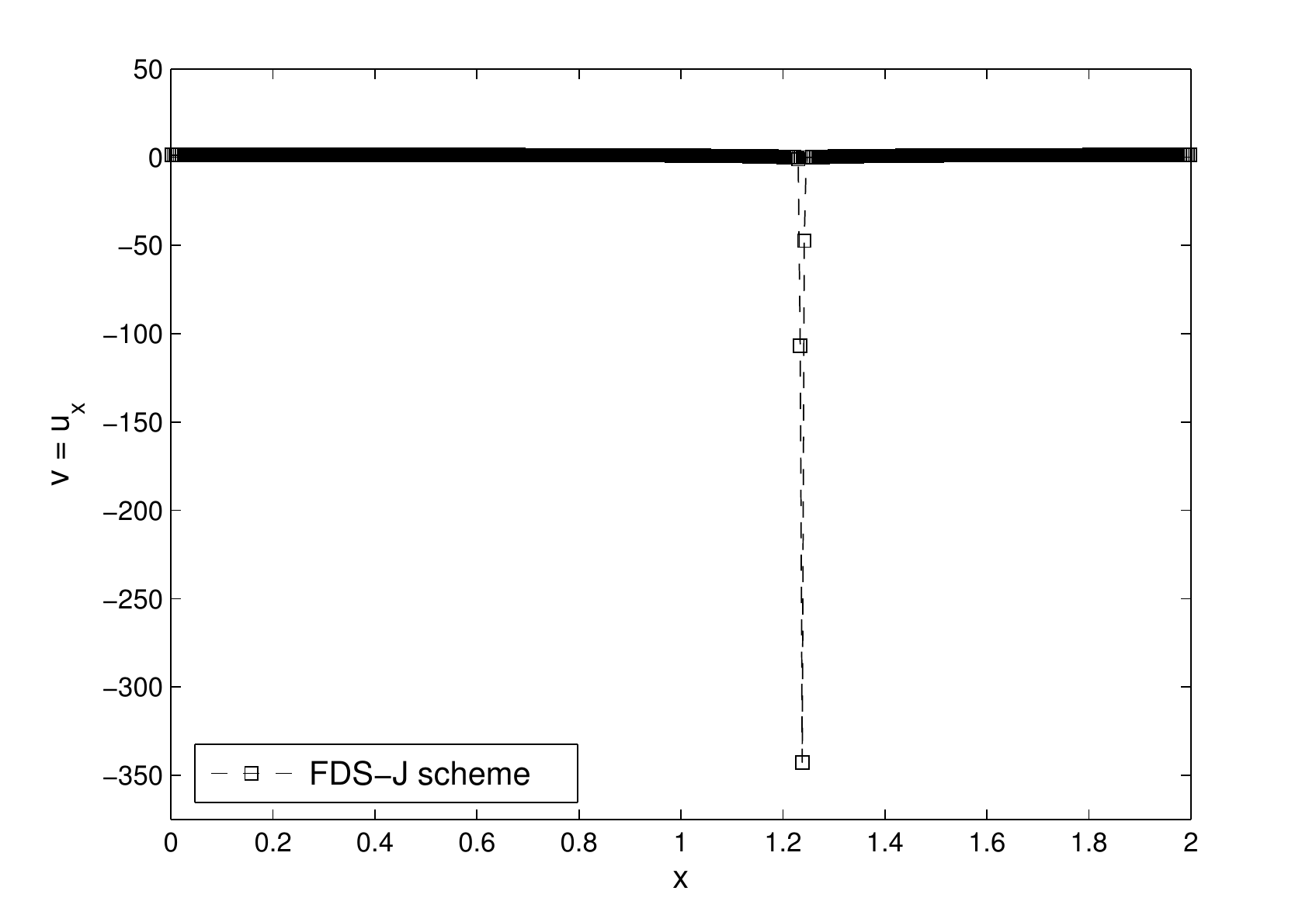}%
\label{Modified Burgers'_delta_shocks}
}%
}%
\caption{ (a) Formation of normal shock in u-variable and (b) represents formation of $\delta-$ shock in v-variable, for Modified Burgers' system.}
\end{figure}

\begin{figure}[!ht]
\centerline{%
\subfigure[]{%
\includegraphics[trim=0 5 35 5, clip, width=0.55\textwidth]{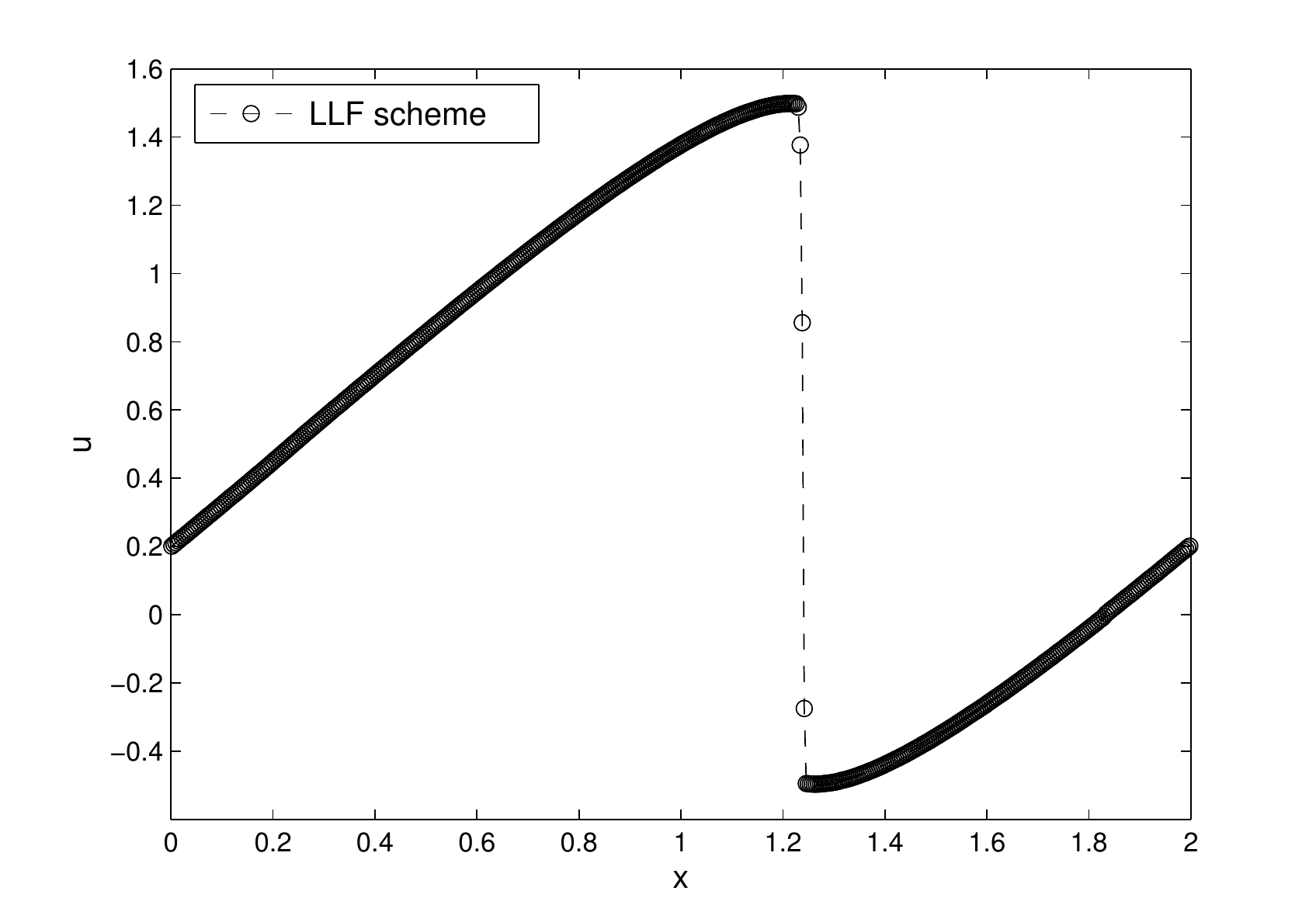}
\label{Modified Burgers'_normal_shock_comp}
}%
\subfigure[]{%
\includegraphics[trim=0 5 35 5, clip, width=0.55\textwidth]{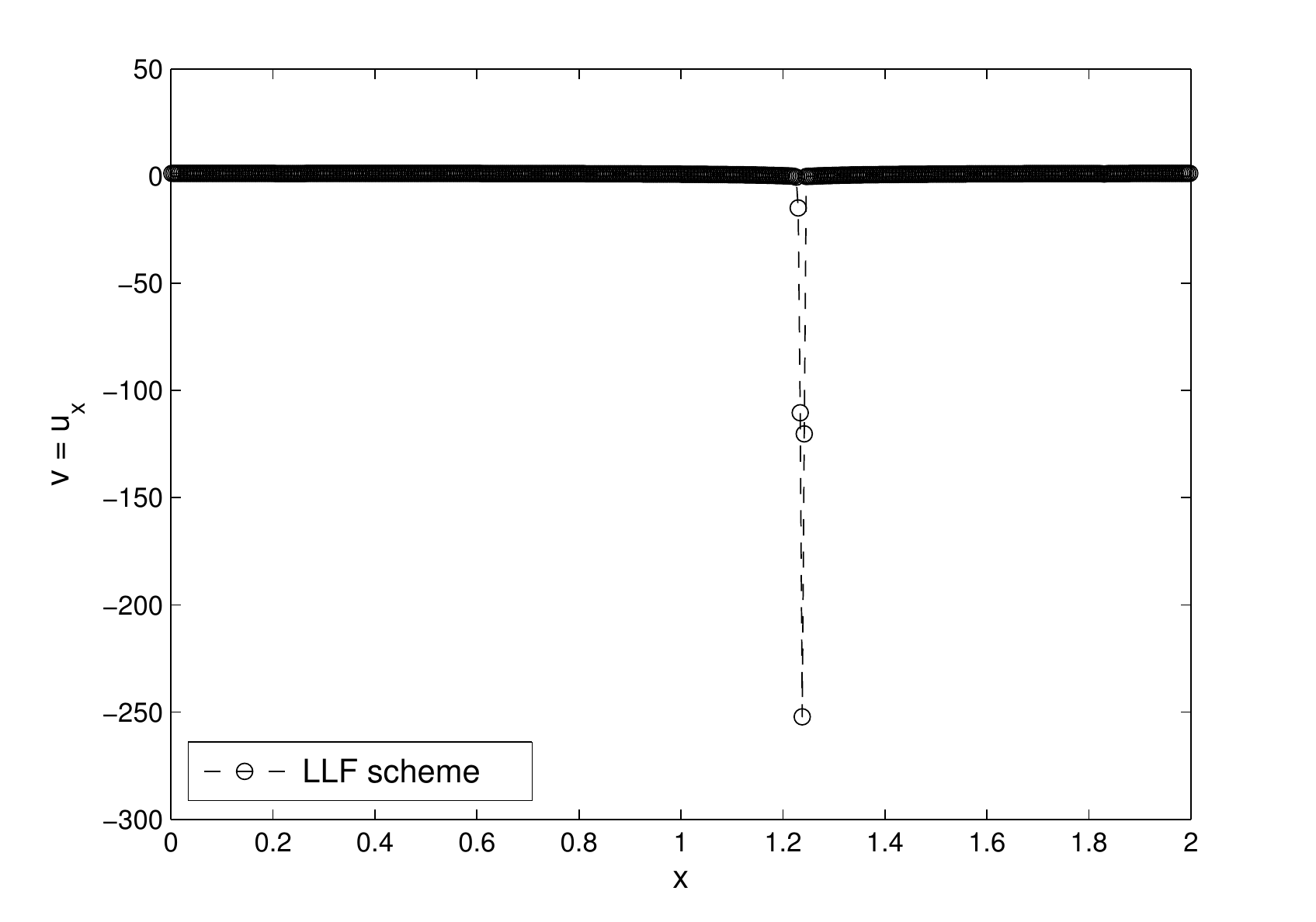}%
\label{Modified Burgers'_delta_shocks_comp}
}%
}%
\caption{ LLF scheme with 500 points (a)  formation of normal shock in u-variable and (b) represents formation of $\delta-$ shock in v-variable, for Modified Burgers' system.}
\end{figure}

\begin{figure}[!ht]
\centerline{%
\subfigure[]{%
\includegraphics[trim=0 5 35 5, clip, width=0.55\textwidth]{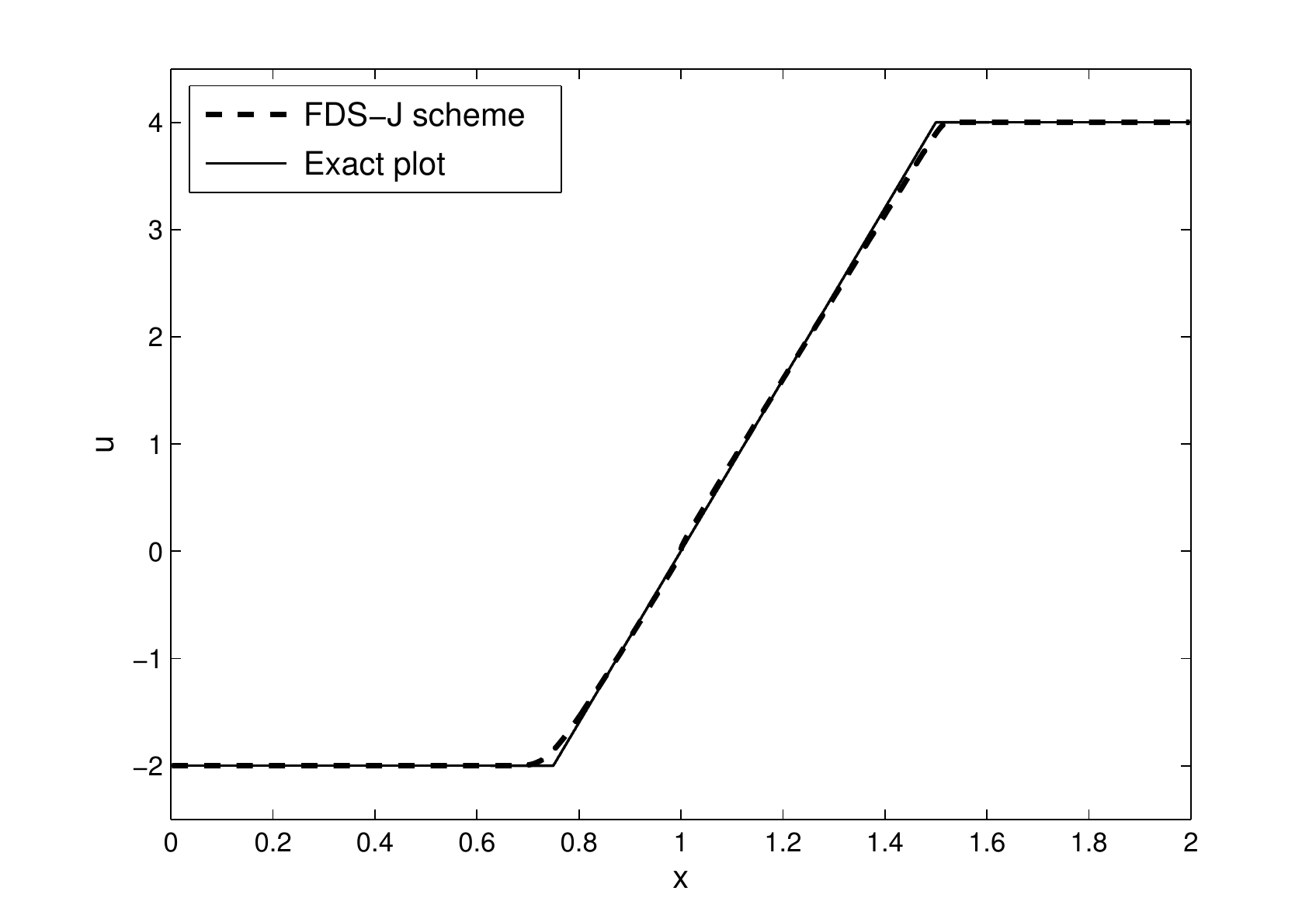}
\label{Modified Burgers'_test2_u}
}%
\subfigure[]{%
\includegraphics[trim=0 5 35 5, clip, width=0.55\textwidth]{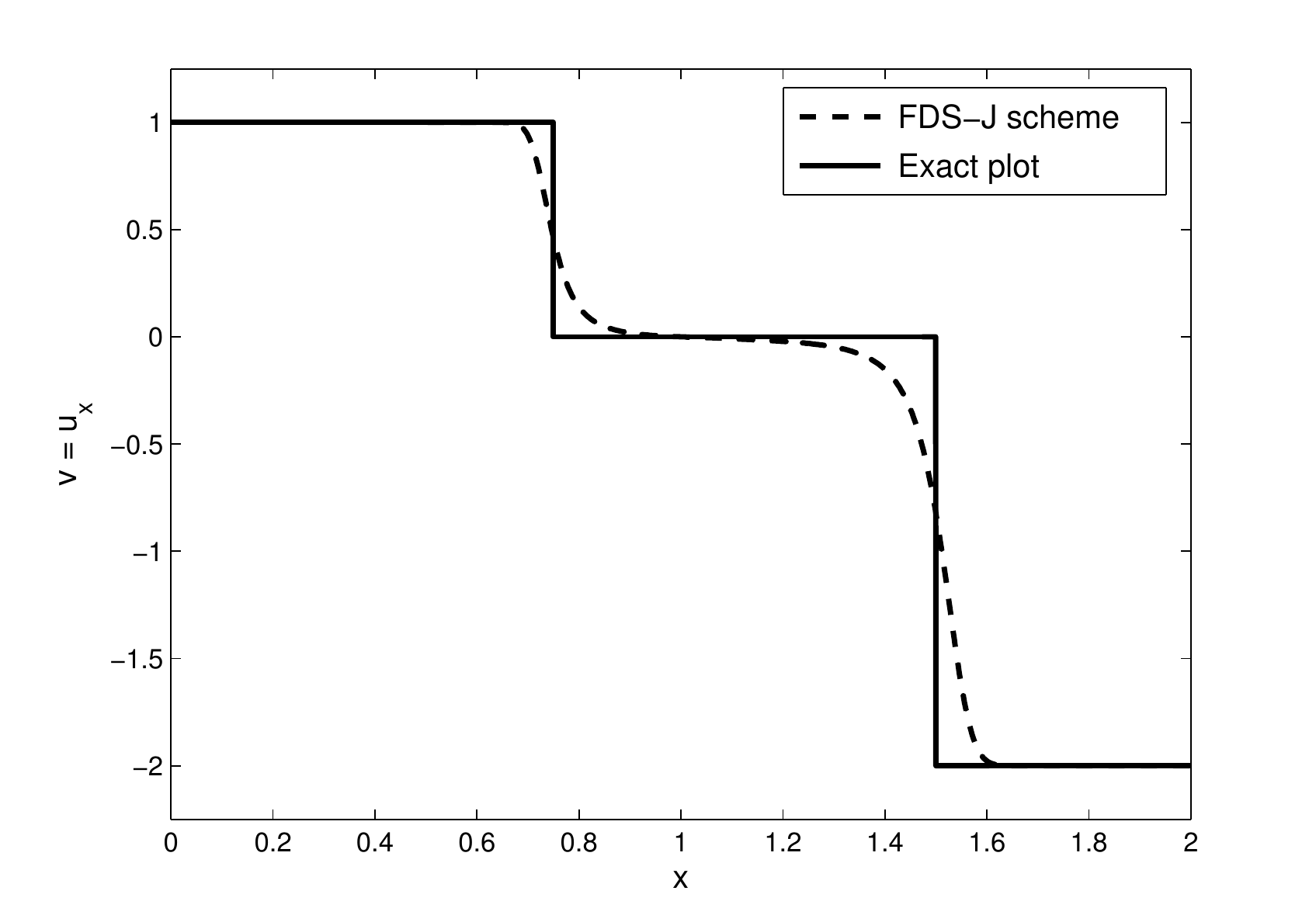}%
\label{Modified Burgers'_test2_v}
}%
}%
\caption{ (a) Sonic point problem results in u-variable and (b) represents numerical results in v-variable for FDS-J scheme with Modified Burgers' system.}
\end{figure}

We considered some numerical test cases from \cite{Capdeville} for the modified Burgers' system.  First test case contains smooth initial conditions which are given as 
\begin{eqnarray*}
 \boldsymbol{U}(x,0) \ = \  \left\{ \begin{array}{l}
 \frac{1}{2} + sin(\pi x) \\ \pi cos(\pi x) \end{array} \right. \forall x \in [0,2] 
\end{eqnarray*} 
with a 2-periodic boundary condition.  Later near time $t =$ $\frac{3}{(2\pi)}$, the given system develops a normal shock and a $\delta-$shock in $u$ and $v$ variables respectively. Theoretically, $v = \pi cos(\pi x)$ may be zero at points $x = \frac{1}{2}, \frac{3}{2}$ but computationally it is not so. Results with FDS-J scheme are given in Figure \ref{Modified Burgers'_normal_shock} and \ref{Modified Burgers'_delta_shocks}. Next we present results with Local Lax-Friedrichs (LLF) method, which is a simple central solver and  are given in Figure \ref{Modified Burgers'_normal_shock_comp}, \ref{Modified Burgers'_delta_shocks_comp}.  Second test case for which initial conditions are defined as $(u_L, v_L) = (-2.0, 1.0)$, $(u_R, v_R) = (4.0, -2.0)$ with $x_{o} = 1.0$ contains a sonic point. Final solutions are obtained at time $t=0.125$ units as given in \cite{Smith_et_al}. Harten's entropy fix is employed to get meaningful solutions and results are given in Figure \ref{Modified Burgers'_test2_u} and \ref{Modified Burgers'_test2_v}.
\section{Further modified Burgers' system}
Shelkovich \cite{Shelkovich} shows existence of $\delta^{\prime}$-shocks in addition to $\delta$-shocks. These shocks occur in  a system which is formed by taking one more derivative of second equation of modified Burgers' system leading to $3\times3$ system. Similarly, Joseph \cite{Joseph} shows existence of $\delta^{\prime\prime}$-shocks in the solution of $4\times4$ system. Let us consider again both equations of modified Burgers' system
\begin{equation}
  u_{t}  \ + \ f_{x}(u) \ = \ 0
 \end{equation}
 and
 \begin{equation}
 v_{t}  \ + \ g_{x}(u) \ = \ 0
 \end{equation}
 On differentiating above equation w.r.t $x$, we get
\begin{equation}
 w_{t}  \ + \ (v^{2} + uw)_{x} \ = \ 0
 \end{equation}
 If we differentiate above equation once more we have
 \begin{equation}
 z_{t}  \ + \ (3vw + uz)_{x} \ = \ 0
 \end{equation}
 In a quasi-linear form above set of four equations can be written as
 \begin{equation}
 \boldsymbol{U}_{t}  \ + \ \boldsymbol{A}\boldsymbol{U}_{x} \ = \ 0
 \end{equation}
 where, $\boldsymbol{U}$ is a $4\times1$ column vector and $\boldsymbol{A}$ is a Jacobian matrix which is given below
 \begin{equation}
  \boldsymbol{A} = \begin{bmatrix}
        \ u  && 0 && 0 && 0  \\[0.3em]
        \ v  &&  u && 0 && 0  \\[0.3em]
        \ w  && 2v && u && 0  \\[0.3em]
        \ z && 3w && 3v && u 
        \end{bmatrix}
 \end{equation}
Eigenvalues corresponding to matrix $\boldsymbol{A}$ are $u,u,u,u$ and for $v\neq0,w\neq0,z\neq0$, matrix $\boldsymbol{A}$ is weakly hyperbolic. Indeed it has only one LI eigenvector $\boldsymbol{e}_{4}$. In this case also we find that there is only one Jordan block of order $4$ as $rank(\boldsymbol{A} - u \boldsymbol{I})^{4}  \ = \ 0 \ = \  rank(\boldsymbol{A} - u \boldsymbol{I})^{5}$. This means for present system, a Jordan chain of order four corresponding to eigenvalue $\lambda = u$ will form, {\em i.e.}, 
\begin{align}
 \begin{split}
  \boldsymbol{A} \boldsymbol{R}_{1}  \ &= \   \lambda \boldsymbol{R}_{1} \\
 \boldsymbol{A}  \boldsymbol{R}_{2}  \ &= \   \lambda \boldsymbol{R}_{2}  \ + \   \boldsymbol{R}_{1} \\
 \boldsymbol{A}  \boldsymbol{R}_{3}  \ &= \  \lambda  \boldsymbol{R}_{3}  \ + \  \boldsymbol{R}_{2} \\
 \boldsymbol{A}  \boldsymbol{R}_{4}  \ &= \  \lambda  \boldsymbol{R}_{4}  \ + \  \boldsymbol{R}_{3}
 \end{split}
\end{align}
where $\boldsymbol{R}_{1} = \boldsymbol{e}_{4}$ and on using $\boldsymbol{R}_{1}$ in the second relation, $\boldsymbol{R}_{2}$ comes out as  $(0,0,\frac{1}{3v},x_{4})^{t}$ with $x_{4}$ as a real constant. Similarly, on using $\boldsymbol{R}_{2}$ in next relation, $\boldsymbol{R}_{3}$ comes out equal to $(0,\frac{1}{6v^{2}},\frac{x_{4}}{3v}-\frac{w}{6v^{3}},y_{4})$, where $x_{4}$ is already defined and $y_{4}$ is another real constant. Finally, last expression gives $\boldsymbol{R}_{4} = (\frac{1}{6v^{3}}, \frac{x_{4}}{6v^{2}}, \frac{y_{4}}{3v} - \frac{z}{18v^{4}} - \frac{w x_{4}}{6v^{3}}, t_{4})^{t}$. Let $\boldsymbol{P}$ denote a matrix with column vectors $[\boldsymbol{R}_{1} | \boldsymbol{R}_{2} | \boldsymbol{R}_{3} | \boldsymbol{R}_{4}]$ and one can check determinant of $\boldsymbol{P}$ is $\frac{1}{108v^{6}} \neq 0$.
\subsection{Formulation of FDS scheme for Further Modified Burgers' System} 
\begin{figure}[!ht] 
\centerline{%
\subfigure[]{%
\includegraphics[trim=0 5 35 5, clip, width=0.55\textwidth]{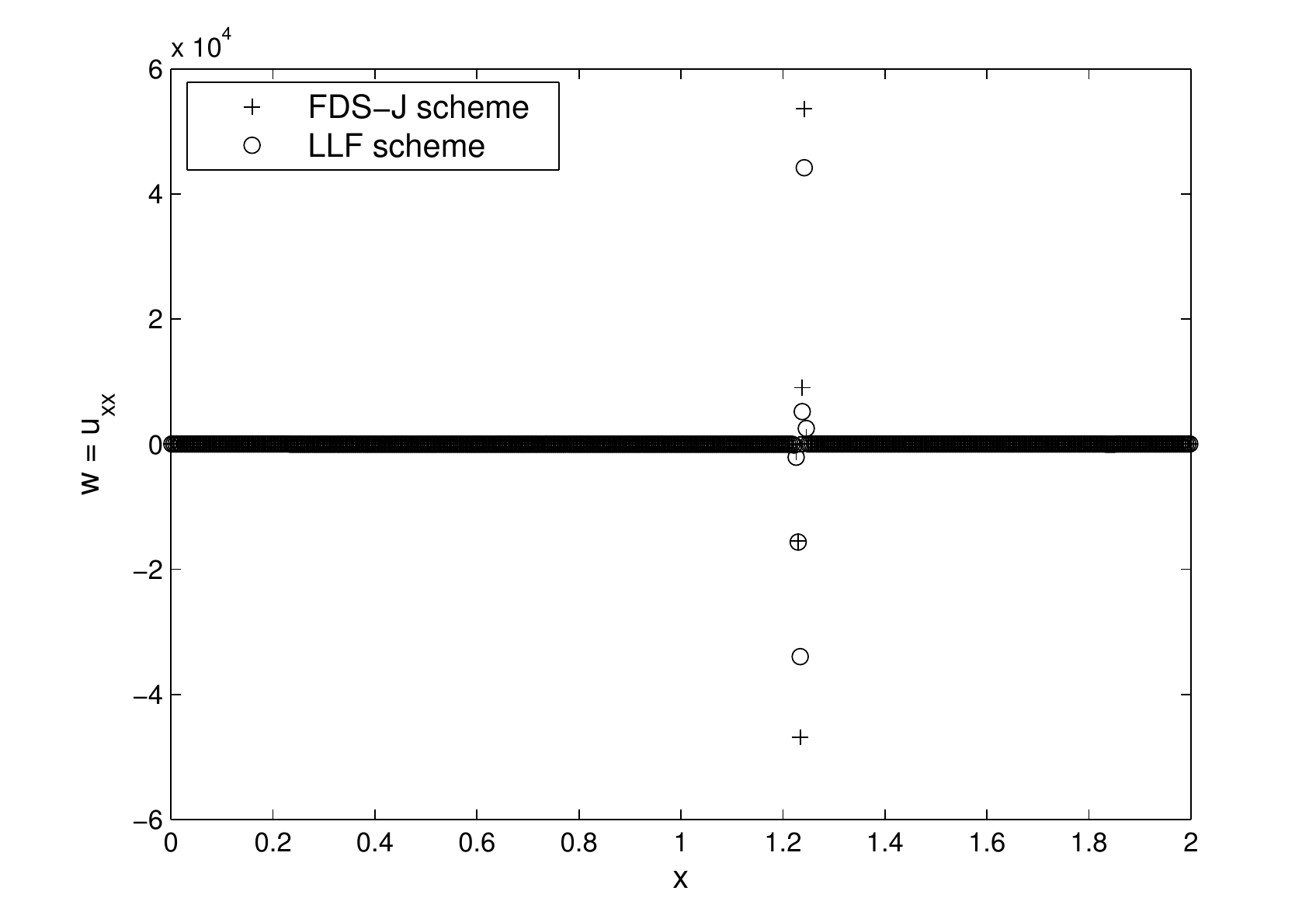}
\label{delta_dash_shock}
}%
\subfigure[]{%
\includegraphics[trim=0 5 35 5, clip, width=0.55\textwidth]{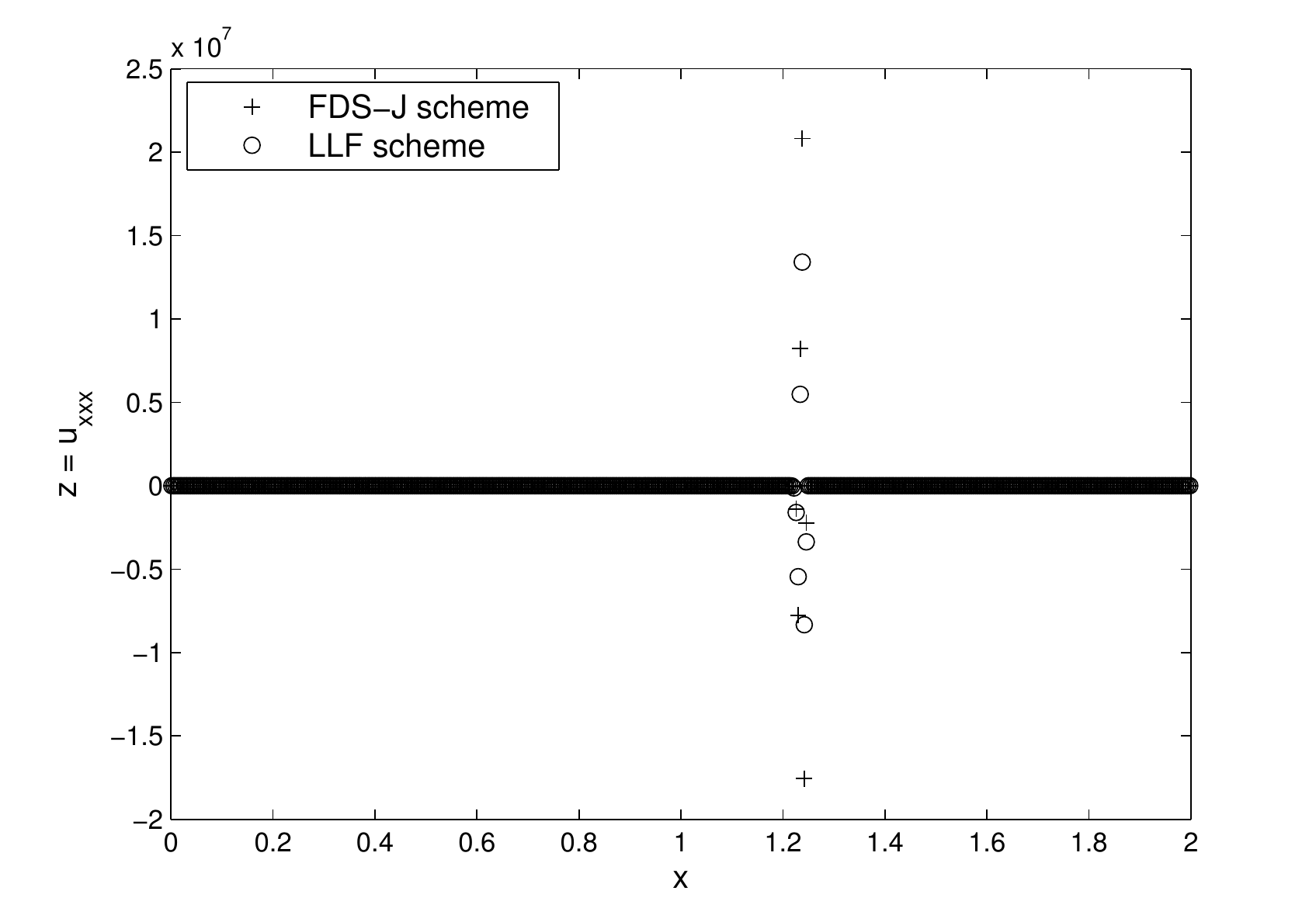}%
\label{delta_double_dash_shock}
}%
}%
\caption{ Comparison of FDS-J scheme with LLF scheme for further modified Burgers' system (a) represents formation of $\delta^{\prime}-$shock in w-variable and (b) represents formation of $\delta^{\prime\prime}-$shock in z-variable.}
\end{figure}

In this case $\bigtriangleup{\boldsymbol{F}}$ is written as,
\begin{equation*} 
\bigtriangleup{\boldsymbol{F}}  \ = \   \bar\alpha_{1} \bar \lambda \boldsymbol{\bar{R}}_{1}  \ + \ \bar\alpha_{2}(\bar \lambda \boldsymbol{\bar{R}}_{2}  \ + \   \boldsymbol{\bar{R}}_{1})  \ + \ \bar\alpha_{3}(\bar\lambda  \boldsymbol{\bar{R}}_{3}  \ + \  \boldsymbol{\bar{R}}_{2}) \ + \ \bar\alpha_{4}(\bar \lambda  \boldsymbol{\bar{R}}_{4}  \ + \  \boldsymbol{\bar{R}}_{3})
\end{equation*} 
After splitting each of the eigenvalues into a positive part and a negative part, $\bigtriangleup{\boldsymbol{F}}^{+}$ and $\bigtriangleup{\boldsymbol{F}}^{-}$ can be written as
\begin{equation}
\bigtriangleup{\boldsymbol{F}}^{+}   \ = \   \bar\alpha_{1} \bar \lambda^{+} \boldsymbol{\bar{R}}_{1}  \ + \ \bar\alpha_{2}(\bar \lambda^{+} \boldsymbol{\bar{R}}_{2}  \ + \   \boldsymbol{\bar{R}}_{1}) \ + \ \bar\alpha_{3}(\bar\lambda^{+}  \boldsymbol{\bar{R}}_{3}  \ + \  \boldsymbol{\bar{R}}_{2}) \ + \ \bar\alpha_{4}(\bar \lambda^{+}  \boldsymbol{\bar{R}}_{4}  \ + \  \boldsymbol{\bar{R}}_{3})
\end{equation}
and
\begin{equation}
\bigtriangleup{\boldsymbol{F}}^{-}   \ = \   \bar\alpha_{1} \bar \lambda^{-} \boldsymbol{\bar{R}}_{1}  \ + \ \bar\alpha_{2}(\bar \lambda^{-} \boldsymbol{\bar{R}}_{2}  \ + \   \boldsymbol{\bar{R}}_{1}) \ + \ \bar\alpha_{3}(\bar\lambda^{-}  \boldsymbol{\bar{R}}_{3}  \ + \  \boldsymbol{\bar{R}}_{2}) \ + \ \bar\alpha_{4}(\bar \lambda^{-}  \boldsymbol{\bar{R}}_{4}  \ + \  \boldsymbol{\bar{R}}_{3})
\end{equation}
\begin{equation}\label{flux_differencing_3}
\Rightarrow \bigtriangleup{\boldsymbol{F}}^{+} - \bigtriangleup{\boldsymbol{F}}^{-}  \ = \   |\bar \lambda|\bigtriangleup{\boldsymbol{U}}  
\end{equation}
In this case $\bigtriangleup{\boldsymbol{U}}$ is defined as,
\begin{equation}
  \bigtriangleup{\boldsymbol{U}} \  =  \begin{bmatrix}
         \ \bigtriangleup{u} \\[0.3em]
         \ \bigtriangleup{v}  \\[0.3em]
         \ \bigtriangleup{w}  \\[0.3em]
         \ \bigtriangleup{z}
       \end{bmatrix}
\end{equation}
and $\bar{\lambda} = \bar{u}$. In order to solve (\ref{flux_differencing_3}) fully, we need to find average value of $u$. In this case also average value of $u$ turns out to be equal to  $\dfrac{u_{L} + u_{R}}{2}$. we take the same test case as considered in the modified Burgers' system with initial smooth conditions
\begin{eqnarray*}
 \boldsymbol{U}(x,0) \ = \ \left\{ \begin{array}{l} \frac{1}{2} + sin(\pi x) \\  \pi cos(\pi x) \\  -\pi^{2} sin(\pi x) \\ -\pi^{3} cos(\pi x)  \end{array} \right. \forall x \in [0,2] 
\end{eqnarray*} 
As already explained at time $t =$ $\frac{3}{(2\pi)}$, the given system develops a normal shock and a $\delta-$shock in $u$ and $v$ variables. Similarly, at same position where normal shock forms, third variable $w$ gives a $\delta^{\prime}$-shock and fourth variable $z$ creates a $\delta^{\prime \prime}$-shock. Results for FDS-J scheme are compared with simple central solver LLF and are given in Figures \ref{delta_dash_shock} and \ref{delta_double_dash_shock}.  
\section{Summary}
In this study, we attempted to develop a Flux Difference Splitting scheme for genuine weakly hyperbolic systems to simulate various shocks including  $\delta$-shocks, $\delta^{\prime}$-shocks and $\delta^{\prime\prime}$-shocks.   Newly constructed FDS-J scheme, developed using Jordan Canonical forms together with an upwind flux difference splitting method, is capable of recognizing these shocks accurately.  For considered weakly hyperbolic systems, there is no direct contribution of generalized eigenvector in the final formulation of the scheme.


\begin{thebibliography}{99}
\bibitem{Bouchut_Jin_&_Li} F. Bouchut, S. Jin and X. Li (2003). Numerical approximations of pressureless and isothermal gas dynamics. SIAM Journal on Numerical Analysis, 41(1), 135-158.
\bibitem{Capdeville} G. Capdeville (2008). Towards a compact high-order method for non-linear hyperbolic systems, II. The Hermite-HLLC scheme. Journal of Computational Physics, 227(22), 9428-9462.
\bibitem{Chen_&_Liu} G. Q. Chen and H. Liu (2003). Formation of $\delta$-shocks and vacuum states in the vanishing pressure limit of solutions to the Euler equations for isentropic fluids. SIAM journal on mathematical analysis, 34(4), 925-938.
\bibitem{Shelkovich} V. G. Danilov and V. M. Shelkovich (2005). Dynamics of propagation and interaction of $\delta$-shock waves in conservation law systems. Journal of Differential Equations, 211(2), 333-381.
\bibitem{Harten_entropy_fix} A. Harten (1984). On a class of high resolution total-variation-stable finite-difference schemes. SIAM Journal on Numerical Analysis, 21(1), 1-23.
\bibitem{Joseph} K. T. Joseph and M. R. Sahoo (2013). Vanishing viscosity approach to a system of conservation laws admitting $\delta^{\prime\prime}$ waves. Communications on Pure \& Applied Analysis, 12(5).
\bibitem{Osher} S.J. Osher and F. Solomon, Upwind difference schemes for hyperbolic systems of conservation laws, Mathematics of Computation, vol. 38, no. 158, pp. 339-374, 1982. 
\bibitem{Roe} P. L. Roe (1981). Approximate Riemann solvers, parameter vectors, and difference schemes. Journal of computational physics, 43(2), 357-372.
\bibitem{LLF} V. V. E. Rusanov (1962). The calculation of the interaction of non-stationary shock waves and obstacles. USSR Computational Mathematics and Mathematical Physics, 1(2), 304-320.
\bibitem{Smith_et_al} T. A. Smith, D. J. Petty and C. Pantano (2016). A Roe-like numerical method for weakly hyperbolic systems of equations in conservation and non-conservation form. Journal of Computational Physics, 316, 117-138.
\end{thebibliography}
\end{document}